\newif\ifspringer
\springerfalse
\ifspringer

    \documentclass{svmult}

    \smartqed
    \usepackage{mathptmx}       % selects Times Roman as basic font
    \usepackage{helvet}         % selects Helvetica as sans-serif font
    \usepackage{courier}        % selects Courier as typewriter font
    \usepackage{type1cm}        % activate if the above 3 fonts are
    \renewcommand{\email}[1]{\emailname: #1} % change the email address font style
    
    \spdefaulttheorem{assumption}{Assumption}{\upshape \bfseries}{\itshape}
    \spdefaulttheorem{algo}{Algorithm}{\upshape \bfseries}{\itshape}

\else

    \documentclass{article}

    \usepackage[a4paper]{geometry}

    \newenvironment{acknowledgement}{\paragraph{Acknowledgement.}}{\par}
    \newcommand{\email}[1]{\texttt{#1}}

    \usepackage{amsthm}
    \theoremstyle{plain}

    \theoremstyle{definition}

    \theoremstyle{remark}

\fi

\usepackage{makeidx}         % allows index generation
\usepackage{graphicx}        % standard LaTeX graphics tool
\usepackage[bottom]{footmisc}% places footnotes at page bottom

\usepackage{latexsym}
\usepackage{amsmath}
\usepackage{amsfonts}
\usepackage{amssymb}
\usepackage{bm}

\usepackage{url}
\usepackage{algorithm}
\usepackage{algorithmic}
\usepackage[misc,geometry]{ifsym}

\newcommand{\bsa}{{\boldsymbol{a}}}

\newcommand{\bst}{{\boldsymbol{t}}}

\newcommand{\bsw}{{\boldsymbol{w}}}
\newcommand{\bsx}{{\boldsymbol{x}}}
\newcommand{\bsy}{{\boldsymbol{y}}}
\newcommand{\bsz}{{\boldsymbol{z}}}

\newcommand{\bsell}{{\boldsymbol{\ell}}}
\newcommand{\bszero}{{\boldsymbol{0}}} % vector of zeros
\newcommand{\bsalpha}{{\boldsymbol{\alpha}}}

\newcommand{\bsgamma}{{\boldsymbol{\gamma}}}

\newcommand{\bsnu}{{\boldsymbol{\nu}}}

\newcommand{\bstau}{{\boldsymbol{\tau}}}

\newcommand{\bsDelta}{{\boldsymbol{\Delta}}}

\newcommand{\rd}{{\mathrm{d}}}

\newcommand{\bbE}{{\mathbb{E}}}

\newcommand{\bbR}{{\mathbb{R}}}

\newcommand{\bbZ}{{\mathbb{Z}}}
\DeclareSymbolFont{bbold}{U}{bbold}{m}{n}
\DeclareSymbolFontAlphabet{\mathbbold}{bbold}

\newcommand{\calA}{{\mathcal{A}}}

\newcommand{\calN}{{\mathcal{N}}}
\newcommand{\calO}{{\mathcal{O}}}

\newcommand{\setu}{{\mathfrak{u}}}
\newcommand{\setv}{{\mathfrak{v}}}

\usepackage{microtype} % good font tricks

\usepackage[colorlinks=true,linkcolor=black,citecolor=black,urlcolor=black]{hyperref}
\usepackage{bookmark}
\makeatletter % to avoid hyperref warnings:
  \providecommand*{\toclevel@author}{999}
  \providecommand*{\toclevel@title}{0}
\makeatother

\begin{document}

\ifspringer

    \title*{Hot New Directions for Quasi-Monte Carlo Research in Step with
    Applications}
    \author{Frances Y.\ Kuo \and Dirk Nuyens}

    \institute{
     Frances Y. Kuo (\Letter)
     \at School of Mathematics and Statistics,
         University of New South Wales, Sydney NSW 2052, Australia \\
     \email{f.kuo@unsw.edu.au}
     \and
     Dirk Nuyens
     \at Department of Computer Science, KU Leuven,
                  Celestijnenlaan 200A, 3001 Leuven, Belgium \\
     \email{dirk.nuyens@cs.kuleuven.be}
    }

    \index{Kuo, Frances Y.}
    \index{Nuyens, Dirk}

\else

    \title{Hot New Directions for Quasi-Monte Carlo Research \\
        in Step with Applications}

    \author{Frances Y.\ Kuo\footnote{
            Frances Y. Kuo (\Letter):
            School of Mathematics and Statistics,
            University of New South Wales, Sydney NSW 2052, Australia,
            \email{f.kuo@unsw.edu.au}
        }
       \and Dirk Nuyens\footnote{
            Dirk Nuyens:
            Department of Computer Science, KU Leuven,
            Celestijnenlaan 200A, 3001 Leuven, Belgium,
            \email{dirk.nuyens@cs.kuleuven.be}
        }
    }

    \date{}

\fi

\maketitle

\abstract{This article provides an overview of some interfaces between the
theory of quasi-Monte Carlo (QMC) methods and applications. We summarize
three QMC theoretical settings: first order QMC methods in the unit cube
$[0,1]^s$ and in $\bbR^s$, and higher order QMC methods in the unit cube.
One important feature is that their error bounds can be independent of the
dimension $s$ under appropriate conditions on the function spaces. Another
important feature is that good parameters for these QMC methods can be
obtained by fast efficient algorithms even when $s$ is large. We outline
three different applications and explain how they can tap into the
different QMC theory. We also discuss three cost saving strategies that
can be combined with QMC in these applications. Many of these recent QMC
theory and methods are developed not in isolation, but in close connection
with applications.}

\section{Introduction}

High dimensional computation is a new frontier in scientific computing,
with applications ranging from financial mathematics such as option
pricing or risk management, to groundwater flow, heat transport, and wave
propagation. A tremendous amount of progress has been made in the past two
decades on the theory and application of \emph{quasi-Monte Carlo}
(\emph{QMC}) \emph{methods} for approximating high dimensional integrals.
See e.g., the classical references \cite{Nie92,SJ94} and the recent books
\cite{DP10,Lem09,LP14}. One key element is the \emph{fast
component-by-component construction} \cite{CKN06,Nuy14,NC06a,NC06b} which
provides parameters for first order or \emph{higher order QMC methods}
\cite{DKS13,DP10} for sufficiently smooth functions. Another key element
is the careful selection of parameters called \emph{weights}
\cite{SWW04,SW98} to ensure that the worst case errors in an appropriately
weighted function space are bounded independently of the dimension. The
dependence on dimension is very much the focus of the study on
\emph{tractability} \cite{NW10} of multivariate problems.

We are particularly keen on the idea that new theory and methods for high
dimensional computation are developed not in isolation, but in close
connection with applications. The theoretical QMC convergence rates depend
on the appropriate pairing between the function space and the class of QMC
methods. Practitioners are free to choose the theoretical setting or
pairing that is most beneficial for their applications, i.e., to achieve
the best possible convergence rates under the weakest assumptions on the
problems. As QMC researchers we take application problems to be our guide
to develop new theory and methods as the needs arise. This article
provides an overview of some interfaces between such theory and
applications.

We begin in Section~\ref{sec:settings} by summarizing three theoretical
settings. The first setting is what we consider to be the standard QMC
setting for integrals formulated over the unit cube. Here the integrand is
assumed to have square-integrable mixed first derivatives, and it is
paired with \emph{randomly shifted lattice rules} \cite{SKJ02b} to achieve
first order convergence. The second setting is for integration over
$\bbR^s$ against a product of univariate densities. Again the integrands
have square-integrable mixed first derivatives and we use randomly shifted
lattice rules to achieve first order convergence. The third setting
returns to the unit cube, but considers integrands with higher order mixed
derivatives and pairs them with \emph{interlaced polynomial lattice rules}
\cite{God15} which achieve higher order convergence. These three settings
are discussed in more detail in \cite{KN16}.

Next in Section~\ref{sec:apps} we outline three applications of QMC
methods: option pricing, GLMM (generalized linear mixed models) maximum
likelihood, PDEs with random coefficients -- all with quite different
characteristics and requiring different strategies to tackle them. We
explain how to match each example application with an appropriate setting
from Section~\ref{sec:settings}. In the option pricing application, see
e.g., \cite{ABG98,GKSW08,GilWat09}, none of the settings is applicable due
to the presence of a \emph{kink}. We discuss the strategy of
\emph{smoothing by preintegration} \cite{GKLS}, which is similar to the
method known as conditional sampling \cite{ACN13}. In the maximum
likelihood application \cite{KDSWW08}, the change of variables plays a
crucial role in a similar way to importance sampling for Monte Carlo
methods. In the PDE application, see e.g., \cite{CD15,GWZ14,KN16,SG11},
the \emph{uniform} and the \emph{lognormal} cases correspond to
integration over the unit cube and $\bbR^s$, respectively, and the two
cases tap into different QMC settings. For the lognormal case we briefly
contrast three ways to generate the random field: \emph{Karhunen--Lo\`eve
expansion}, \emph{circulant embedding}
\cite{DN97,GKNSS11,GKNSS-paper1,GKNSS-paper2}, and \emph{$H$-matrix
technique} \cite{FKS,Hac15}.

Then in Section~\ref{sec:savings} we discuss three different cost saving
strategies that can be applied to all of the above applications. First,
\emph{multi-level methods} \cite{Gil15} restructure the required
computation as a telescoping sum and tackle different levels separately to
improve the overall cost versus error balance, while more general
\emph{multi-index methods} \cite{HNT16} allow different criteria to be
considered simultaneously in a multi-index telescoping sum. Second, the
\emph{multivariate decomposition methods} \cite{KNPSW,KSWW10a,Was13} work
in a similar way by making an explicit decomposition of the underlying
function into functions of only subsets of the variables \cite{KSWW10b}.
The third strategy is \emph{fast QMC matrix-vector multiplication} which
carries out the required computation for multiple QMC samples at the same
time using an FFT \cite{DKLS15}.

\begin{figure} \label{fig:overview}
\begin{center}
    \ifspringer
        \includegraphics[width=1\textwidth]{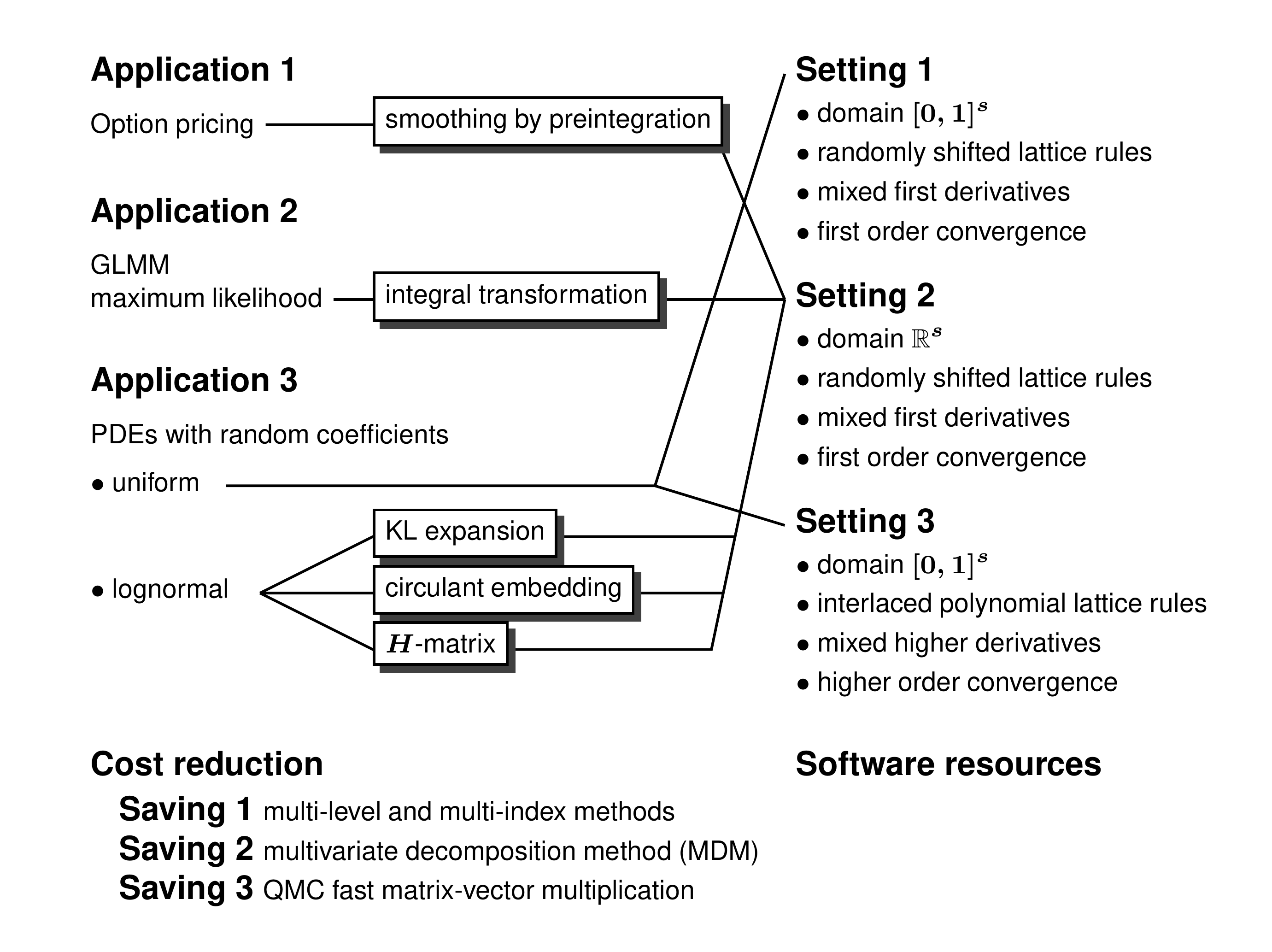}
    \else
        \includegraphics[width=0.8\textwidth]{overview}
    \fi
\end{center}
\caption{The connection between different components of this article.}
\end{figure}

We provide pointers to some software resources in
Section~\ref{sec:software} and conclude the article in
Section~\ref{sec:summary} with a summary and an outlook to future work. An
overview of the various components of this article is given in
Figure~\ref{fig:overview}.

\section{Three Settings} \label{sec:settings}

Here we describe three theoretical function space settings paired with
appropriate QMC methods. These three setting are also covered in
\cite{KN16}. Of course these three pairs are not the only possible
combinations. We selected them due to our preference for constructive QMC
methods that achieve the best possible convergence rates, with the implied
constant independent of dimension, under the weakest possible assumptions
on the integrands.

\subsection{Setting 1: Standard QMC for the Unit Cube} \label{sec:setting1}

For $f$ a real-valued function defined over the $s$-dimensional unit cube
$[0,1]^s$, with $s$ finite and fixed, we consider the integral
\begin{align} \label{eq:int1}
  I(f)
  \,=\, \int_{[0,1]^s} f(\bsy)\, \rd\bsy\,.
\end{align}

\subsubsection*{Weighted Sobolev spaces}

We assume in this standard setting that the integrand $f$ belongs to a
\emph{weighted Sobolev space of smoothness one} in the unit cube
$[0,1]^s$. Here we focus on the \emph{unanchored} variant in which the
norm is defined by, see also \cite{SWW04},
\begin{align} \label{eq:norm1}
  \|f\|_\bsgamma
  \,=\,
  \Bigg[
  \sum_{\setu\subseteq\{1:s\}}
  \frac{1}{\gamma_\setu}
  \int_{[0,1]^{|\setu|}}
  \bigg(\int_{[0,1]^{s-|\setu|}}
  \frac{\partial^{|\setu|}f}{\partial \bsy_\setu}(\bsy)
  \,\rd\bsy_{\{1:s\}\setminus\setu}
  \bigg)^2
  \rd\bsy_\setu
  \Bigg]^{1/2},
\end{align}
where $\{1:s\}$ is a shorthand notation for the set of indices
$\{1,2,\ldots,s\}$, $(\partial^{|\setu|}f)/(\partial \bsy_\setu)$ denotes
the mixed first derivative of $f$ with respect to the ``active'' variables
$\bsy_\setu = (y_j)_{j\in\setu}$, while $\bsy_{\{1:s\}\setminus\setu} =
(y_j)_{j\in\{1:s\}\setminus\setu}$ denotes the ``inactive'' variables.

There is a weight parameter $\gamma_\setu\ge 0$ associated with each
subset of variables~$\bsy_\setu$ to model their relative importance. We
denote the weights collectively by $\bsgamma$. Special forms of weights
have been considered in the literature. \emph{POD weights} (\emph{product
and order dependent weights}), arisen for the first time in \cite{KSS12},
take the form
\[
  \gamma_\setu \,=\, \Gamma_{|\setu|}\,\prod_{j\in\setu} \Upsilon_j\,,
\]
which is specified by two sequences $\Gamma_0= \Gamma_1=1,
\Gamma_2,\Gamma_3,\ldots\ge 0$ and $\Upsilon_1\ge\Upsilon_2\ge \cdots
> 0$. Here the factor $\Gamma_{|\setu|}$ is
said to be order dependent because it is determined solely by the
cardinality of $\setu$ and not the precise indices in $\setu$. The
dependence of the weight $\gamma_\setu$ on the indices $j\in\setu$ is
controlled by the product of terms $\Upsilon_j$. Each term $\Upsilon_j$ in
the sequence corresponds to one coordinate direction; the sequence being
non-increasing indicates that successive coordinate directions become less
important. Taking all $\Gamma_{|\setu|}=1$ or all $\Upsilon_j = 1$
corresponds to the weights known as \emph{product weights} or \emph{order
dependent weights} in the literature \cite{SWW04,SW98}.

\subsubsection*{Randomly shifted lattice rules}

We pair the weighted Sobolev space with \emph{randomly shifted lattice
rules}; the complete theory can be found in \cite{DKS13}. Randomly shifted
lattice rules approximate the integral \eqref{eq:int1} by
\begin{align} \label{eq:lat}
  Q(f)
  \,=\, \frac{1}{n} \sum_{i=1}^n f(\bst_i),\qquad
  \bst_i \,=\, \left\{\frac{i\,\bsz}{n} + \bsDelta\right\}\,,
\end{align}
where $\bsz\in \bbZ^s$ is known as the \emph{generating vector},
$\bsDelta$ is a \emph{random shift} drawn from the uniform distribution
over $[0,1]^s$, and the braces indicate that we take the fractional parts
of each component in a vector.

A randomly shifted lattice rule provides an unbiased estimate of the
integral, i.e., $\bbE[Q(f)] = I(f)$, where the expectation is taken with
respect to the random shift $\bsDelta$. Its quality is determined by the
choice of the generating vector $\bsz$. By analyzing the quantity known as
\emph{shift-averaged worst case error}, it is known that good generating
vectors can be obtained using a \emph{CBC construction}
(\emph{component-by-component construction}), determining the components
of $\bsz$ one at a time sequentially, to achieve nearly $\calO(n^{-1})$
convergence rate which is optimal in the weighted Sobolev space of
smoothness one, and the implied constant in the big $\calO$ bound can be
independent of $s$ under appropriate conditions on the weights $\bsgamma$.

More precisely, if $n$ is a power of $2$ then we know that the CBC
construction yields the root-mean-square error bound, for all $\lambda\in
(1/2,1]$,
\begin{align} \label{eq:error1}
  \sqrt{\bbE\left[ |I(f) - Q(f)|^2 \right]}
  &\,\le\,
  \Bigg(
  \frac{2}{n}
  \sum_{\emptyset\ne\setu\subseteq\{1:s\}} \gamma_\setu^\lambda\,
  [\vartheta(\lambda)]^{|\setu|}
  \Bigg)^{1/(2\lambda)}
  \,\|f\|_\bsgamma\,,
\end{align}
where $\vartheta(\lambda) := 2\zeta(2\lambda)/(2\pi^2)^\lambda$, with
$\zeta(a) := \sum_{k=1}^\infty k^{-a}$ denoting the Riemann zeta function.
A similar result holds for general $n$. The best rate of convergence
clearly comes from choosing $\lambda$ close to $1/2$, but the advantage is
offset by the fact that $\zeta(2\lambda)\to\infty$ as $\lambda\to
(1/2)_+$.

\subsubsection*{Choosing the weights}

To apply this abstract theory to a given practical integrand $f$, we need
to first obtain an estimate of the norm $\|f\|_\bsgamma$. Remember that at
this stage we do not yet know how to choose the weights $\gamma_\setu$.
Assuming that bounds on the mixed first derivatives in \eqref{eq:norm1}
can be obtained so that
\begin{equation} \label{eq:norm-bound}
  \|f\|_\bsgamma \,\le\,
  \left(\sum_{\setu\subseteq\{1:s\}} \frac{B_\setu}{\gamma_\setu}\right)^{1/2},
\end{equation}
we can substitute \eqref{eq:norm-bound} into \eqref{eq:error1} and then,
with $\lambda$ fixed but unspecified at this point and $A_\setu =
[\vartheta(\lambda)]^{|\setu|}$, we choose the weights $\gamma_\setu$ to
minimizing the product
\[
 C_\bsgamma \,:=\, \left(\sum_{\setu\subseteq\{1:s\}} \gamma_\setu^\lambda\,A_\setu\right)^{1/(2\lambda)}
 \left(\sum_{\setu\subseteq\{1:s\}} \frac{B_\setu}{\gamma_\setu}\right)^{1/2}.
\]
Elementary calculus leads us to conclude that we should take
\begin{equation} \label{eq:weights}
  \gamma_\setu \,:=\, \left(\frac{B_\setu}{A_\setu}\right)^{1/(1+\lambda)},
\end{equation}
which yields
\[
  C_\bsgamma \,=\,
  \Bigg(
  \sum_{\setu\subseteq\{1:s\}} A_\setu^{1/(1+\lambda)} B_\setu^{\lambda/(1+\lambda)}
  \Bigg)^{(1+\lambda)/(2\lambda)}.
\]
We then specify a value of $\lambda$, as close to $1/2$ as possible, to
ensure that $C_\bsgamma$ can be bounded independently of $s$. This in turn
determines the theoretical convergence rate which is
$\calO(n^{-1/(2\lambda)})$.

The chosen weights $\gamma_\setu$ are then fed into the CBC construction
to produce generating vectors for randomly shifted lattice rules that
achieve the desired theoretical error bound for this integrand. This
strategy for determining weights was first considered in \cite{KSS12}.

\emph{Fast} CBC constructions (using FFT) can produce generating vectors
for an $n$-point rule in $s$ dimensions in $\calO(s\,n\log n)$ operations
in the case of product weights \cite{Nuy14}, and in $\calO(s\,n\log n +
s^2\,n)$ operations in the case of POD weights \cite{KSS11}. Note that
these are considered to be pre-computation costs. The actual cost for
generating the points on the fly is $\calO(s\,n)$ operations, no worse
than Monte Carlo simulations. Strategies to improve on the computational
cost of approximating the integral are discussed in
Section~\ref{sec:savings}.

The CBC construction yields a lattice rule which is extensible in
dimension~$s$. We can also construct \emph{lattice sequences} which are
extensible or embedded in the number of points~$n$, at the expense of
increasing the implied constant in the error bound
\cite{CKN06,DPW08,HHLL00,HN03}.

\subsection{Setting 2: QMC Integration over $\bbR^s$} \label{sec:setting2}

QMC approximation to an integral which is formulated over the Euclidean
space $\bbR^s$ can be obtained by first mapping the integral to the unit
cube as follows:
\begin{align} \label{eq:int-Rs}
  I(f) \,=\, \int_{\bbR^s} f(\bsy)\,\prod_{j=1}^s \phi(y_j)\,\rd\bsy
  &\,=\, \int_{[0,1]^s} f(\Phi^{\mbox{-}1}(\bsw))\,\rd\bsw \\
  &\,\approx\, \frac{1}{n} \sum_{i=1}^n f(\Phi^{\mbox{-}1}(\bst_i)) \,=\, Q(f) \,. \nonumber
\end{align}
(With a slight abuse of notation we have reused $I(f)$ and $Q(f)$ from the
previous subsection for integration over $\bbR^s$ in this subsection.)
Here $\phi$ can be \emph{any general univariate probability density
function}, and $\Phi^{\mbox{-}1}$ denotes the component-wise application
of the inverse of the cumulative distribution function corresponding to
$\phi$. Note that in many practical applications we need to first apply
some clever transformation to convert the integral into the above form;
some examples are discussed in Section~\ref{sec:apps}. The transformed
integrand $f\circ \Phi^{\mbox{-}1}$ arising from practical applications
typically does not belong to the Sobolev space defined over the unit cube
due to the integrand being unbounded near the boundary of the cube, or
because the mixed derivatives of the transformed integrand do not exist or
are unbounded. Thus the theory in the preceding subsection generally does
not apply in practice. Some theory for QMC on singular integrands is given
in \cite{Owe06}.

We summarize here a special weighted space setting in $\bbR^s$ for which
randomly shifted lattice rules have been shown to achieve nearly the
optimal convergence rate of order one \cite{KSWWat10,NK14}. The norm in
this setting is given by
\begin{align} \label{eq:norm2}
  \|f\|_\bsgamma
  &\,=\,
  \Bigg[
  \sum_{\setu\subseteq\{1:s\}} \frac{1}{\gamma_\setu}
  \int_{\bbR^{|\setu|}}
  \bigg(
  \int_{\bbR^{s-|\setu|}}
  \frac{\partial^{|\setu|} f}{\partial \bsy_\setu}(\bsy)
  \bigg(\prod_{j\in\{1:s\}\setminus\setu} \phi(y_j)\bigg)
  \,\rd\bsy_{\{1:s\}\setminus\setu}
  \bigg)^2 \nonumber\\
  &\qquad\qquad\qquad\qquad\qquad\qquad\qquad\qquad\qquad\;\times
  \bigg(\prod_{j\in\setu} \varpi_j^2(y_j)\bigg)
  \,\rd\bsy_\setu \Bigg]^{1/2} \,.
\end{align}
Comparing~\eqref{eq:norm2} with~\eqref{eq:norm1}, apart from the
difference that the integrals are now over the unbounded domain, there is
a probability density function $\phi$ as well as additional \emph{weight
functions} $\varpi_j$ which can be chosen to moderate the tail behavior of
the mixed derivatives of $f$.

The convergence results for the CBC construction of randomly shifted
lattice rules in this general setting depend on the choices of $\phi$ and
$\varpi_j$. For $n$ a power of $2$, the root-mean-square error bound takes
the form, for all $\lambda\in (1/(2r),1]$,
\[
  \sqrt{\bbE\left[ |I(f) - Q(f)|^2 \right]}
  \,\le\,
  \Bigg(
  \frac2n
  \sum_{\emptyset\ne \setu\subseteq\{1:s\}}
  \gamma_\setu^\lambda\, \prod_{j\in\setu} \vartheta_j(\lambda)\Bigg)^{1/(2\lambda)}\,
   \,
   \|f\|_{\bsgamma}\,,
\]
with $r$ (appearing in the applicable lower bound on $\lambda$) and
$\vartheta_j(\lambda)$ depending on $\phi$ and $\varpi_j$, see
\cite[Theorem~8]{NK14}. Some special cases have been analyzed:
\begin{itemize}
\item %
See \cite[Theorem~15]{GKNSSS15} or \cite[Theorem~5.2]{KN16} for
$\phi(y) = \phi_{\rm nor}(y) = \exp(-y^2/2)/\sqrt{2\pi}$ being the
standard normal density and $\varpi_j^2(y_j) = \exp (-2\,\alpha_j\,
|y_j|)$ with $\alpha_j > 0$.
\smallskip
\item %
See \cite[Theorem~5.3]{KN16} for $\phi = \phi_{\rm nor}$ and
$\varpi_j^2(y_j) = \exp (-\alpha\, y_j^2)$ with $\alpha < 1/2$.
\smallskip
\item %
See \cite[Theorem~2]{SKS13} for $\phi$ being a logistic, normal, or
Student density and $\varpi_j = 1$.
\end{itemize}

To apply this abstract theory to a practical integral over $\bbR^s$, it is
important to realize that the choice of $\phi$ can be tuned as part of the
process of transformation to express the integral in the form
\eqref{eq:int-Rs}. (This point will become clearer when we describe the
maximum likelihood application in Subsection~\ref{sec:app2}.) Then the
choice of weight functions $\varpi_j$ arises as part of the process to
obtain bounds on the norm of~$f$, as in \eqref{eq:norm-bound}. (This point
will become clearer when we describe the PDE application in
Subsection~\ref{sec:app3}.) Finally we can choose the weights
$\gamma_\setu$ as in \eqref{eq:weights} but now with $A_\setu =
\prod_{j\in\setu} \vartheta_j(\lambda)$ for the appropriate
$\vartheta_j(\lambda)$ corresponding to the choice of $\phi$
and~$\varpi_j$. The choice of density $\phi$, weight functions $\varpi_j$,
and weight parameters $\gamma_\setu$ then enter the CBC construction to
obtain the generating vector of good randomly shifted lattice rules that
can achieve the theoretical error bound for this integrand.

In practice, it may well be that the weights $\gamma_\setu$ obtained in
this way are not sensible because we were working with theoretical upper
bounds on the error that may be too pessimistic. It may already be so in
the standard setting of the previous subsection, but is more pronounced in
the setting for $\bbR^s$ due to the additional complication associated
with the presence of $\phi$ and $\varpi_j$.

\subsection{Setting 3: Smooth Integrands in the Unit Cube}
\label{sec:setting3}

Now we return to the integration problem over the unit cube
\eqref{eq:int1} and outline a weighted function space setting from
\cite{DKLNS14} for smooth integrands of order $\alpha$. The norm is given
by
\begin{align}\label{eq:norm3}
 \|f\|_{\bsgamma}
 \,=\,
 \sup_{\setu\subseteq\{1:s\}}
 \sup_{\bsy_\setv \in [0,1]^{|\setv|}}
 \frac{1}{\gamma_\setu}
 \sum_{\setv\subseteq\setu} \,
 \sum_{\bstau_{\setu\setminus\setv} \in \{1:\alpha\}^{|\setu\setminus\setv|}}
 \bigg|\int_{[0,1]^{s-|\setv|}}
 (\partial^{(\bsalpha_\setv,\bstau_{\setu\setminus\setv},\bszero)} f)(\bsy) \,
 \rd \bsy_{\{1:s\} \setminus\setv}
 \bigg|\,
  .
\end{align}
Here $(\bsalpha_\setv,\bstau_{\setu\setminus\setv},\bszero)$ denotes a
multi-index $\bsnu$ with $\nu_j = \alpha$ for $j\in\setv$, $\nu_j =
\tau_j$ for $j\in\setu\setminus\setv$, and $\nu_j = 0$ for $j\notin\setu$.
We denote the $\bsnu$-th partial derivative of $f$ by $\partial^\bsnu f
\,=\,
(\partial^{|\bsnu|}f)/(\partial^{\nu_1}_{y_1}\partial^{\nu_2}_{y_2}\cdots\partial^{\nu_s}_{y_s})$.

This function space setting can be paired with \emph{interlaced polynomial
lattice rules} \cite{God15,GD15} to achieve higher order convergence rates
in the unit cube. A \emph{polynomial lattice rule} \cite{Nie92} is similar
to a lattice rule (see \eqref{eq:lat} without the random shift
$\bsDelta$), but instead of a generating vector of integers we have a
generating vector of polynomials, and thus the regular multiplication and
division are replaced by their polynomial equivalents. We omit the
technical details here. An interlaced polynomial lattice rule with $n =
2^m$ points in $s$ dimensions with interlacing factor $\alpha$ is obtained
by taking a polynomial lattice rule in $\alpha\,s$ dimensions and then
interlacing the bits from every successive $\alpha$ dimensions to yield
one dimension. More explicitly, for $\alpha = 3$, given three coordinates
$x = (0.x_1 x_2 \ldots x_m)_2$, $y = (0.y_1 y_2 \ldots y_m)_2$ and $z =
(0.z_1 z_2 \ldots z_m)_2$ we interlace their bits to obtain $w = (0.x_1
y_1 z_1 x_2 y_2 z_2 \ldots x_m y_m z_m)_2$.

An interlaced polynomial lattice rule with interlacing factor $\alpha\ge
2$, with irreducible modulus polynomial of degree $m$, and with $n=2^m$
points in $s$ dimensions, can be constructed by a CBC algorithm such that,
for all $\lambda \in (1/\alpha,1]$,
\begin{align*}
  |I(f) - Q(f)|
  \,\le\,
  \left(\frac{2}{n} \sum_{\emptyset\ne\setu\subseteq{\{1:s\}}}
  \gamma_\setu^\lambda\, [\vartheta_{\alpha}(\lambda)]^{|\setu|}\right)^{1/\lambda}\,
  \|f\|_{\bsgamma}\,,
\end{align*}
where $\vartheta_{\alpha}(\lambda) :=
  2^{\alpha\lambda(\alpha-1)/2}
  ([1+1/(2^{\alpha\lambda}-2)]^\alpha-1)$.
This result can be found in \cite[Theorem~5.4]{KN16}, which was obtained
from minor adjustments of \cite[Theorem~3.10]{DKLNS14}.

Given a practical integrand $f$, if we can estimate the corresponding
integrals involving the mixed derivatives in \eqref{eq:norm3}, then we can
choose the weights $\gamma_\setu$ so that every term in the supremum is
bounded by a constant, say, $c$. This strategy in \cite{DKLNS14} led to a
new form of weights called \emph{SPOD weights} (\emph{smoothness-driven
product and order dependent weights}); they take the form
\[
 \gamma_\setu
 \,=\,
 \sum_{\bsnu_\setu \in \{1:\alpha\}^{|\setu|}}
 \Gamma_{|\bsnu_\setu|}  \prod_{j\in\setu} \Upsilon_j(\nu_j)\,.
\]
If the weights $\bsgamma$ are SPOD weights, then the fast CBC construction
of the generating vector has cost $\calO(\alpha\,s\, n\log n +
\alpha^2\,s^2 n)$ operations. If the weights $\bsgamma$ are product
weights, then the CBC algorithm has cost $\calO(\alpha\,s\, n\log n)$
operations.

\section{Three Applications} \label{sec:apps}

Integrals over $\bbR^s$ often arise from practical applications in the
form of multivariate expected values
\begin{equation} \label{eq:QoI}
  \bbE_\rho [q] \,=\, \int_{\bbR^s} q(\bsy) \, \rho(\bsy) \, \rd \bsy\,,
\end{equation}
where $q$ is some quantity of interest which depends on a vector $\bsy =
(y_1,\ldots,y_s)$ of parameters or variables in $s$ dimensions, and $\rho$
is some multivariate probability density function describing the
distribution of $\bsy$, \emph{not necessarily a product of univariate
functions} as we assumed in~\eqref{eq:int-Rs}, and so we need to make an
appropriate transformation to apply our theory. Below we discuss three
motivating applications with quite different characteristics, and we will
explain how to make use of the different settings in
Section~\ref{sec:settings}.

\subsection{Application 1: Option Pricing} \label{sec:app1}

Following the Black--Scholes model, integrals arising from option pricing
problems take the general form \eqref{eq:QoI}, with
\[
  q(\bsy) = \max(\mu(\bsy),0)
  \quad\mbox{and}\quad
  \rho(\bsy) = \frac{\exp(-\frac{1}{2}\bsy^{\tt T}\Sigma^{-1}\bsy)}{\sqrt{(2\pi)^s \det(\Sigma)}},
\]
where the variables $\bsy = (y_1,\ldots,y_s)^{\tt T}$ correspond to a
discretization of the underlying Brownian motion over a time interval
$[0,T]$, and the covariance matrix has entries $\Sigma_{ij} =
(T/s)\min(i,j)$. For example, in the case of an \emph{arithmetic average
Asian call option} \cite{ABG98,GKSW08,GilWat09}, the payoff function $q$
depends on the smooth function $\mu(\bsy) = (1/s)\sum_{j=1}^s
S_{t_j}(\bsy) - K$ which is the difference between the average of the
asset prices $S_{t_j}$ at the discrete times and the strike price $K$.

The widely accepted strategy to rewrite these option pricing integrals
from the form \eqref{eq:QoI} to the form \eqref{eq:int-Rs} with product
densities is to take a factorization $\Sigma = AA^{\tt T}$ and apply a
change of variables $\bsy = A\bsy'$. This yields an integral of the form
\eqref{eq:int-Rs} with
\[
  f(\bsy') = q(A\bsy')
  \quad\mbox{and}\quad
  \phi = \phi_{\rm nor}.
\]
The choice of factorization therefore determines the function $f$. For
example, $A$ can be obtained through Cholesky factorization (commonly
known as the \emph{standard construction}; in this case it is equivalent
to generating the Brownian motions sequentially in time), through
\emph{Brownian bridge construction} \cite{CMO97}, or eigenvalue
decomposition sometimes called the \emph{principal component construction}
\cite{ABG98}. Note that in practice these factorizations are not carried
out explicitly due to the special form of the covariance matrix. In fact,
they can be computed in $\calO(s)$, $\calO(s)$ and $\calO(s\,\log s)$
operations, respectively \cite{GKSW08}.

The success of QMC for option pricing cannot be explained by most existing
theory due to the \emph{kink} in the integrand induced by the maximum
function. However, for some factorizations it is shown in \cite{GKS13}
that all \emph{ANOVA terms} of $f$ are smooth, with the exception of the
highest order term. This hints at a \emph{smoothing by preintegration}
strategy, where a coordinate with some required property is chosen, say
$y_k$, and we integrate out this one variable (either exactly or
numerically with high precision) to obtain a function in $s-1$ variables
\[
  P_k(f) \,:=\, \int_{-\infty}^\infty f(\bsy)\,\phi_{\rm nor}(y_k)\,\rd y_k\,.
\]
Under the right conditions (e.g., integrating with respect to $y_1$ in the
case of the principal components construction), this new function is
smooth and belongs to the function space setting of
Subsection~\ref{sec:setting2} (with one less variable) \cite{GKLS}. This
strategy is related to the method known as \emph{conditional sampling}
\cite{ACN13}.

\subsection{Application 2: Maximum Likelihood} \label{sec:app2}

Another source of integrands which motivated recent developments in the
function space setting of Subsection~\ref{sec:setting2} is a class of
generalized linear mixed models (GLMM) in statistics, as examined in
\cite{KDSWW08,KSWWat10,SKS13}. A specific example of the Poisson
likelihood time series model considered in these papers involves an
integral of the form \eqref{eq:QoI}, with
\[
  q(\bsy) = \prod_{j=1}^s \frac{\exp(\tau_j (\beta+y_j) - e^{\beta + y_j})}{\tau_j !}
  \quad\mbox{and}\quad
  \rho(\bsy) = \frac{\exp(-\frac{1}{2}\bsy^{\tt T}\Sigma^{-1}\bsy)}{\sqrt{(2\pi)^s\det(\Sigma)}}.
\]
Here $\beta\in\bbR$ is a model parameter, $\tau_1,\ldots,\tau_s\in
\{0,1,\ldots\}$ are the count data, and $\Sigma$ is a Toeplitz covariance
matrix with $\Sigma_{ij}=\sigma^2\varkappa^{|i-j|}/(1-\varkappa^2)$, where
$\sigma^2$ is the variance and $\varkappa\in(-1,1)$ is the autoregression
coefficient.

An obvious way to rewrite this integral in the form \eqref{eq:int-Rs} with
product densities is to factorize $\Sigma$ as discussed in the previous
subsection for the option pricing applications, but this would yield a
very spiky function~$f$. Instead, the strategy developed in \cite{KDSWW08}
recenters and rescales the exponent $T(\bsy)$ of the product
$q(\bsy)\rho(\bsy) =: \exp(T(\bsy))$ as follows:
\begin{enumerate}
\item %
Find the unique stationary point $\bsy^*$ satisfying $\nabla
T(\bsy^*)=0$.

\item %
Determine the matrix $\Sigma^*=(-\nabla^2T(\bsy^*))^{-1}$ which
describes the convexity of $T$ around the stationary point.

\item %
Factorise $\Sigma^*=A^*{A^*}^{{\tt T}}$.

\item %
Apply a change of variables $\bsy=A^*\bsy'+\bsy^*$.

\item %
Multiply and divide the resulting integrand by the product
$\prod_{j=1}^s \phi(y_j')$ where $\phi$ is any univariate density (not
necessarily the normal density).
\end{enumerate}
These steps then yield an integral of the form \eqref{eq:int-Rs} with
\[
  f(\bsy') \,=\, \frac{c\,\exp (T(A^*\bsy'+\bsy^*))}{\prod_{j=1}^s \phi(y_j')}
\]
for some scaling constant $c>0$. Note that the choice of $A^*$ and $\phi$
determines $f$.

The paper \cite{SKS13} provides careful estimates of the norm of the
resulting integrand~$f$ in the setting of Subsection~\ref{sec:setting2}
corresponding to three different choices of density~$\phi$, with the
weight functions taken as $\varpi_j = 1$, and gives the formula for the
weight parameters $\gamma_\setu$ that minimize the overall error bound.

These GLMM problems are extremely challenging not only for QMC but also in
general the tools are still lacking. There is still lots of room to
develop new QMC methods and theory for these problems.

\subsection{Application 3: PDEs with Random Coefficients} \label{sec:app3}

Our third application is motivated by fluid flow through a porous medium,
typically modelled using Darcy's Law, with random coefficients. A popular
toy problem is the elliptic PDE with a random coefficient
\cite{CD15,GWZ14,SG11}
\[
 -\nabla \cdot (a(\bsx, \omega)\,\nabla u(\bsx, \omega))
 = \kappa(\bsx)\quad \text{for $\bsx\in D\subset\bbR^d$ and almost all $\omega\in\Omega$, }
\]
with $d\in\{1,2,3\}$, subject to homogeneous Dirichlet boundary
conditions. The coefficient $a(\bsx, \omega)$ is assumed to be a random
field over the spatial domain $D$ (e.g., representing the permeability of
a porous material over $D$), and $\Omega$ is the probability space. The
goal is to compute the expected values $\bbE [ G(u) ]$ of some bounded
linear functional $G$ of the solution $u$ over $\Omega$.

For practical reasons it is often assumed that $a(\bsx,\omega)$ is a
\emph{lognormal random field}, that is, $a(\bsx, \omega) = \exp( Z(\bsx,
\omega) )$, where $Z(\bsx, \omega)$ is a Gaussian random field with a
prescribed mean and covariance function. This is known as the
\emph{lognormal case}. However, researchers often analyze a simpler model
known as the \emph{uniform case}.

\subsubsection*{The uniform case}

In the uniform case, we consider the parametric PDE
\begin{equation} \label{eq:PDE-par}
 -\nabla \cdot (a(\bsx, \bsy)\,\nabla u(\bsx, \bsy))
 = \kappa(\bsx)\quad \text{for $\bsx\in D\subset\bbR^d$},
\end{equation}
together with
\begin{equation} \label{eq:uniform}
  a(\bsx, \bsy) \,=\, a_0(\bsx) + \sum_{j\ge1} y_j \psi_j(\bsx)\,,
\end{equation}
where the parameters $y_j$ are independently and uniformly distributed on
the interval $[-\frac{1}{2},\frac{1}{2}]$, and we assume that $0<
a_{\min}\le a(\bsx,\bsy)\le a_{\max} < \infty$ for all $\bsx$ and $\bsy$.

A (single-level) strategy for approximating $\bbE[G(u)]$ is as follows:
\begin{enumerate}
\item Truncate the infinite sum in \eqref{eq:uniform} to $s$ terms.
\item Solve the PDE using finite element methods with meshwidth $h$.
\item Approximate the resulting $s$-dimensional integral using QMC
    with $n$ points.
\end{enumerate}
So the error is a sum of truncation error, discretization error, and
quadrature error.

For the QMC quadrature error in Step~3, we have the integral
\eqref{eq:int1} with
\[
 f(\bsy) \,=\, G(u^s_h(\cdot,\bsy - \tfrac{\boldsymbol{1}}{\boldsymbol{2}})),
\]
where $u^s_h$ denotes the finite element solution of the truncated
problem, and the subtraction by $\tfrac{\boldsymbol{1}}{\boldsymbol{2}}$
takes care of the translation from the usual unit cube $[0,1]^s$ to
$[-\frac{1}{2},\frac{1}{2}]^s$. By differentiating the PDE
\eqref{eq:PDE-par}, we can obtain bounds on the mixed derivatives of the
PDE solution with respect to $\bsy$. This leads to bounds on the
norm~\eqref{eq:norm1} of the integrand $f$ and so we can apply the
theoretical setting of Subsection~\ref{sec:setting1} to obtain up to first
order convergence for QMC. Under appropriate assumptions and with first
order finite elements, we can prove that the total error for the above
3-step strategy is of order \cite{KSS12}
\[
  \calO(s^{-2(1/p_0-1)} + h^{2} + n^{-\min(1/p_0-1/2,1-\delta)})\,,\quad \delta\in (0,\tfrac{1}{2})\,,
\]
where $p_0\in (0,1)$ should be as small as possible while satisfying
$\sum_{j\ge 1} \|\psi_j\|_{L_\infty}^{p_0} < \infty$. This part is
presented as a step-by-step tutorial in the article \cite{KN-tutorial}
from this volume.

The bounds on the derivatives of the PDE with respect to $\bsy$ also allow
us to obtain bounds on the norm \eqref{eq:norm3} and so we can also apply
the theoretical setting of Subsection~\ref{sec:setting3} to obtain higher
order convergence \cite{DKLNS14}. Specifically, the
$\calO(n^{-\min(1/p_0-1/2,1-\delta)})$ term can be improved to
$\calO(n^{-1/p_0})$. Also the $\calO(h^2)$ term can be improved by using
higher order finite elements. See \cite{KN16,KN-tutorial} for more
details.

\subsubsection*{The lognormal case with Karhunen--Lo\`{e}ve expansion}

In the lognormal case, we have the same parametric PDE \eqref{eq:PDE-par},
but now we use the \emph{Karhunen--Lo\`{e}ve expansion} (\emph{KL
expansion}) of the Gaussian random field (in the exponent) to write
\[
 a(\bsx,\bsy) \,=\, a_0(\bsx)\,\exp\left(\sum_{j\ge1} y_j\,\sqrt{\mu_j} \, \xi_j(\bsx)\right)\,,
\]
where $a_0(\bsx)>0$, the $\mu_j$ are real, positive and non-increasing in
$j$, the $\xi_j$ are orthonormal in $L_2(D)$, and the parameters
$y_j\in\bbR$ are standard $\calN(0,1)$ random variables. Truncating the
infinite series in $a(\bsx,\bsy)$ to $s$ terms and solving the PDE with a
finite element method as in the uniform case, we have now an integral of
the form \eqref{eq:int-Rs} with
\[
   f(\bsy) \,=\, G(u^s_h(\cdot,\bsy))
   \quad\mbox{and}\quad \phi = \phi_{\rm nor}\,.
\]
One crucial step in the analysis of \cite{GKNSSS15} is to choose suitable
weight functions $\varpi_j$ so that the function $f$ has a finite and
indeed small norm~\eqref{eq:norm2}, so that the theoretical setting of
Subsection~\ref{sec:setting2} can be applied. Again see
\cite{KN16,KN-tutorial} for more details.

In this lognormal case with KL expansion (and also the uniform case), the
cost per sample of the random field is $\calO(s\,M)$ operations, where $M$
is the number of finite element nodes. This dominates the cost in
evaluating the integrand function under the assumption that assembling the
stiffness matrix to solve the PDE (which depends on the random field) is
higher than the cost of the PDE solve which is $\calO(M\,\log M)$. When
$s$ is large the cost of sampling the random field can be prohibitive, and
this is why the following alternative strategies emerged.

\subsubsection*{The lognormal case with circulant embedding}

Since we have a Gaussian random field we can actually sample the random
field exactly on any set of $M$ spatial points. This leads to an integral
of the form \eqref{eq:QoI} with (assuming the field has zero mean)
\[
  q(\bsy) \,=\, G(u^s_h(\cdot,\bsy))
  \quad\mbox{and}\quad
  \rho(\bsy) = \frac{\exp(-\frac{1}{2}\bsy^{\tt T} R^{-1}\bsy)}{\sqrt{(2\pi)^s \det(R)}},
\]
where $R$ is an $M\times M$ covariance matrix, and initially we have $s=
M$. (Note the subtle abuse of notation that the second argument in
$u^s_h(\bsx,\cdot)$ has a different meaning to the KL case, in the sense
that there the covariance is already built in.) This integral can be
transformed into the form \eqref{eq:int-Rs} with a factorization $R =
AA^{\tt T}$ and a change of variables $\bsy = A\bsy'$, as in the option
pricing example, to obtain
\[
  f(\bsy') \,=\, G(u^s_h(\cdot,A\bsy'))
  \quad\mbox{and}\quad
  \phi = \phi_{\rm nor}\,.
\]
The advantage of this discrete formulation is that there is no error
arising from the truncation of the KL expansion. However, the direct
factorization and matrix-vector multiplication require $\calO(M^3)$
operations which can be too costly when $M$ is large.

The idea of \emph{circulant embedding}
\cite{DN97,GKNSS11,GKNSS-paper1,GKNSS-paper2} is to sample the random
field on a regular grid and to embed the covariance matrix of these points
into a larger $s\times s$ matrix which is \emph{nested block circulant
with circulant blocks}, so that FFTs can be used to reduce the per sample
cost to $\calO(s\,\log s)$ operations. Values of the random field at the
finite element quadrature nodes can be obtained by interpolation. Note
that this turns the problem into an even higher dimensional integral, and
we can have $s\gg M$. For this strategy to work we need to use regular
spatial grid points to sample the field and a stationary covariance
function (i.e., the covariance depends only on the relative distance
between points). An additional difficulty is to ensure positive
definiteness of the extended matrix; this is studied
in~\cite{GKNSS-paper1}.

\subsubsection*{The lognormal case with $H$-matrix technique}

Another approach for the discrete matrix formulation of the lognormal case
is to first approximate $R$ by an $H$-matrix \cite{Hac15} and make use of
$H$-matrix techniques to compute the matrix-vector multiplication with the
square-root of this $H$-matrix at essentially linear cost $\calO(M)$. Two
iterative methods have been proposed in \cite{FKS} to achieve this (one is
based on a variant of the \emph{Lanczos iteration} and the other on the
\emph{Schultz iteration}), with full theoretical justification for the
error incurred in the $H$-matrix approximation. An advantage of this
approach over circulant embedding is that it does not require the spatial
grid to be regular nor that the covariance be stationary.

\subsubsection*{Other developments}

A different QMC analysis for the lognormal case has been considered in
\cite{HPS17}. QMC for holomorphic equations was considered in
\cite{DLS16}, and for Baysesian inversion in \cite{DGLS,SST17}. Recently
there is also QMC analysis developed for the situation where the functions
in the expansion of $a(\bsx,\bsy)$ have local support, see
\cite{GHS,HerSch,Kaz}.

\section{Three Cost Saving Strategies} \label{sec:savings}

In this section we discuss the basic ideas of three different kinds of
cost saving strategies that can be applied to QMC methods, without going
into details. Actually, the circulant embedding and $H$-matrix technique
discussed in the previous section can also be considered as cost saving
strategies. These strategies are not mutually exclusive, and it may be
possible to mix and match them to benefit from compound savings.

\subsection{Saving 1: Multi-level and Multi-index} \label{sec:sav1}

The \emph{multi-level} idea \cite{Gil15} is easy to explain in the context
of numerical integration. Suppose that there is a sequence
$(f_\ell)_{\ell\ge 0}$ of approximations to an integrand $f$, with
increasing accuracy and cost as $\ell$ increases, such that we have the
telescoping sum
\[
 f \,=\, \sum_{\ell=0}^\infty (f_\ell - f_{\ell-1}), \quad
 f_{-1}:=0.
\]
For example, the different $f_\ell$ could correspond to different number
of time steps in option pricing, different number of mesh points in a
finite element solve for PDE, different number of terms in a KL expansion,
or a combination of aspects. A multi-level method for approximating the
integral of $f$ is
\[
 A_{\rm ML}(f) \,=\, \sum_{\ell=0}^L Q_\ell(f_\ell - f_{\ell-1}),
\]
where the parameter $L$ determines the number of \emph{levels}, and for
each level we apply a different quadrature rule $Q_\ell$ to the difference
$f_\ell - f_{\ell-1}$.

The integration error (in this simple description with deterministic
quadrature rules) satisfies
\[
 |I(f) - A_{\rm ML}(f)| \;\le\;
 \underbrace{|I(f) - I_L(f_L)|\vphantom{\sum_{\ell=0}^\infty }}_{\le\, \varepsilon/2}
 \;+\; \underbrace{\sum_{\ell=0}^L |(I_\ell-Q_\ell)(f_\ell - f_{\ell-1})|}_{\le\, \varepsilon/2}.
\]
For a given error threshold $\varepsilon>0$, the idea (as indicated by the
underbraces) is that we choose $L$ to ensure that the first term (the
truncation error) on the right-hand side is $\le\varepsilon/2$, and we
specify parameters for the quadrature rules $Q_\ell$ so that the second
term (the quadrature error) is also $\le\varepsilon/2$. The latter can be
achieved with a Lagrange multiplier argument to minimize cost subject to
the given error threshold. Our hope is that the successive differences
$f_\ell-f_{\ell-1}$ will become smaller with increasing $\ell$ and
therefore we would require less quadrature points for the more costly
higher levels.

The \emph{multi-index} idea \cite{HNT16} generalizes this from a scalar
level index $\ell$ to a vector index $\bsell$ so that we can vary a number
of different aspects (e.g., spatial/temporal discretization)
simultaneously and independently of each other. It makes use of the sparse
grid concept so that the overall cost does not blow up with respect to the
dimensionality of $\bsell$, i.e., the number of different aspects being
considered. A simple example is that we use different finite element
meshwidths for different spatial coordinates. This is equivalent to
applying sparse finite element methods within a multilevel algorithm, see
the article \cite{GilKS} in this volume.

Multi-level and multi-index extensions of QMC methods for the applications
from Section~\ref{sec:apps} include e.g.,
\cite{DKLS16,GilWat09,KSSSU,KSS15,RNV}.

\subsection{Saving 2: Multivariate Decomposition Method} \label{sec:sav2}

In the context of numerical integration, the \emph{multivariate
decomposition method} (\emph{MDM}) \cite{GilWas17,KNPSW,KSWW10a,Was13}
makes use of a decomposition of the integrand $f$ of the form
\[
  f \,=\, \sum_{|\setu|<\infty} f_\setu\,,
\]
where the sum is over all finite subsets $\setu\subset\{1,2,\ldots\}$ and
each function $f_\setu$ depends only on the integration variables with
indices in the set $\setu$. Then MDM takes the form
\[
 A_{\rm MDM}(f) \,=\, \sum_{\setu\in\calA} Q_\setu(f_\setu)\,
\]
where $\calA$ is known as the \emph{active set} of subsets of indices, and
for each $\setu$ in the active set we apply a different quadrature rule
$Q_\setu$ to $f_\setu$.

Analogously to the multi-level idea, the error of MDM satisfies
\[
  |I(f) - A_{\rm MDM}(f)| \;\le\;
 \underbrace{\sum_{\setu\notin\calA}|I_\setu(f_\setu)|}_{\le\, \varepsilon/2}
 \;+\; \underbrace{\sum_{\setu\in\calA} |(I_\setu-Q_\setu)(f_\setu)|}_{\le\, \varepsilon/2}\,,
\]
where we choose the active set $\calA$ to ensure that the truncation error
is $\le\varepsilon/2$, and we use a Lagrange multiplier argument to
specify parameters for the quadrature rules so that the quadrature error
is also $\le\varepsilon/2$. Our hope is that, although the cardinality of
the active set $\calA$ might be huge (e.g., tens of thousands), the
cardinality of the individual subsets $\setu\in\calA$ might be relatively
small (e.g., at most $8$ or $10$), and therefore we transfer the problem
into that of solving a large number of low dimensional integrals.

There are many important considerations in the implementation of MDM
\cite{GKNW}. First, we need to decide on how to decompose the integrand
$f$ so that values of the functions $f_\setu$ can be computed. One obvious
choice is known as the \emph{anchored decomposition} which can be computed
via the explicit formula \cite{KSWW10b}
\begin{equation} \label{eq:anchored}
 f_\setu(\bsy_\setu) \,=\,
 \sum_{\setv\subseteq\setu} (-1)^{|\setu|-|\setv|}
 f(\bsy_\setv;\bsa),
\end{equation}
where $\bsa$ is an \emph{anchor} and $(\bsy_\setv;\bsa)$ denotes a vector
obtained from $\bsy$ by replacing the component $y_j$ with the
corresponding component $a_j$ when the index $j$ does not belong to the
subset $\setv$. (This is similar to the well-known ANOVA decomposition
which, however, involves integrals that cannot be computed in practice.)
Second, we need to specify and construct the active set $\calA$ and have
an efficient data structure to store the sets for later traversing. Third,
we need to explore nestedness or embedding in the quadrature rules, taking
into account the sum in \eqref{eq:anchored} and develop efficient ways to
reuse function evaluations.

\subsection{Saving 3: Fast QMC Matrix-vector Multiplication} \label{sec:sav3}

There is a certain structure in some QMC methods that can allow for fast
matrix-vector multiplication using FFT. This structure has been exploited
in the fast CBC construction of lattice rules and polynomial lattice rules
\cite{Nuy14}. We now explain how this same structure can also be used in
more general circumstances \cite{DKLS15}.

For notational convenience, we denote all QMC points $\bst_i$ as row
vectors in this subsection. Given an arbitrary matrix $A$, suppose we want
to
\[
  \mbox{compute}\quad
  \bsy_i\,A
  \quad\mbox{for all}\quad i=1,\ldots,n\,,
\]
with the row vectors $\bsy_i = \chi(\bst_i)$, where $\chi$ denotes an
arbitrary univariate function that is applied to every component of the
QMC point $\bst_i$. Typically we have $\bst_n \equiv \bst_0 = \bszero$ so
we can leave it out. Consider for simplicity the case $n$ is prime and
suppose we can write
\[
  Y\,:=\,
 \begin{bmatrix} \bsy_1 \\\vdots \\ \bsy_{n-1}\end{bmatrix} \,=\, C\, P\,
\]
where $C$ is an $(n-1)\times(n-1)$ circulant matrix, while $P$ is a matrix
containing a single $1$ in each column and $0$ everywhere else. Then we
can compute $Y\bsa$ in $\calO(n\log n)$ operations for any column $\bsa$
of $A$.

The desired factorization $Y = CP$ is possible if we have deterministic
lattice points or deterministic polynomial lattice points, and if we apply
the inverse cumulative distribution function mapping or tent transform
\cite{CKNS16,DNP14,Hic02}. However, it does \emph{not} work with random
shifting, scrambling \cite{Owe97}, or interlacing. This strategy can be
used to generate normally distributed points with a general covariance
matrix (no need for stationarity as in circulant embedding), solving PDEs
with uniform random coefficients, or solving PDEs with lognormal random
coefficients involving finite element quadratures.

\section{Software Resources} \label{sec:software}

We provide some software resources for the practical application of QMC
methods:

\begin{itemize}
\item \emph{The Magic Point Shop}: a collection of QMC point
    generators and
    generating vectors. \\
\url{https://people.cs.kuleuven.be/~dirk.nuyens/qmc-generators/}
\smallskip

\item \emph{Fast component-by-component constructions}: a collection
    of software routines for fast CBC constructions of generating
    vectors.  \\
    \url{https://people.cs.kuleuven.be/~dirk.nuyens/fast-cbc/}
\smallskip

\item \emph{QMC4PDE}: accompanying software package for the survey
    \cite{KN16} on using QMC methods for parametrized
    PDE problems. \\
    \url{https://people.cs.kuleuven.be/~dirk.nuyens/qmc4pde/}
\smallskip

\item \emph{A practical guide to QMC methods}:
   a non-technical introduction of QMC methods with software demos. \\
    \url{https://people.cs.kuleuven.be/~dirk.nuyens/taiwan/}
\end{itemize}

\section{Summary and Outlook} \label{sec:summary}

In this article we summarized three QMC theoretical settings: randomly
shifted lattice rules achieving first order convergence in the unit cube
and in $\bbR^s$, and interlaced polynomial lattice rules achieving higher
order convergence in the unit cube. One important feature is that the
error bound can be independent of the dimension under appropriate
conditions on the weights. Another important feature is that these QMC
methods can be constructed by fast CBC algorithms.

We outlined three different applications and explained how they can be
pre-processed to make use of the different theory. We also discussed three
cost saving strategies that can be combined with QMC in these
applications.

This paper is not meant to be a comprehensive survey on QMC methods. There
are of course many other significant developments on QMC methods and their
applications. For example, we did not discuss \emph{tent transformation}
(also known as the baker's transform), which can yield second order
convergence for randomly shifted rules or first order convergence for
deterministic lattice rules \cite{CKNS16,DNP14,Hic02}. We also did not
discuss \emph{scrambling} \cite{Owe97}, which is a well-known
randomization method that can potentially improve the convergence rates by
an extra half order.

For the future we would like to see QMC in new territories, to tackle a
significantly wider range of more realistic problems. Some emerging new
application areas of QMC include e.g., Bayesian inversion
\cite{DGLS,SST17}, stochastic wave propagation \cite{GH15,GanKS}, quantum
field theory \cite{AGHJLV16,JLGM14}, and neutron transport
\cite{GGKSS,GPS}.

Looking ahead into future QMC developments, what would be on the top of
our wish list? We would very much like to have a \emph{``Setting 4'' where
we have QMC methods that achieve higher order convergence in $\bbR^s$,
with error bounds that are independent of $s$, and for which fast
constructions are possible.} This open problem has seen some partial
solutions \cite{DILP,NN} but there is more to be done!

%%%%%%%%%%%%%%%%%%%%%%%%%%%%%%%%%%%%%%%%%%%%%%%%%%%%%%%%%%%%%%%%%%%%%%%%%%%%%%%%%%%%%%%%%%%
%%% The acknowledgements
\begin{acknowledgement}
The authors acknowledge the financial supports from the Australian
Research Council (FT130100655 and DP150101770) and the KU Leuven research
fund (OT:3E130287 and C3:3E150478).
\end{acknowledgement}

%%%%%%%%%%%%%%%%%%%%%%%%%%%%%%%%%%%%%%%%%%%%%%%%%%%%%%%%%%%%%%%%%%%%%%%%%%%%%%%%%%%%%%%%%%%
%%% The bibliography
%
% BibTeX users please use
\bibliographystyle{spmpsci}
\bibliography{mybibfile}

\begin{thebibliography}{999}

\bibitem{ACN13}
  %{N.~Achtsis, R.~Cools, and D.~Nuyens},
  {Achtsis, N., Cools, R., Nuyens, D.:}
  {Conditional sampling for barrier option pricing under the Heston model},
  in: Monte Carlo and Quasi-Monte Carlo Methods 2012 (J.~Dick, F.~Y.~Kuo, G.~Peters, and I.~H.~Sloan, eds),
  Springer, Berlin, pp.~253--269 (2013)

\bibitem{ABG98}
  %{P.~Acworth, M.~Broadie, and P.~Glasserman},
  {Acworth, P., Broadie, M., Glasserman, P.:}
  {A comparison of some Monte Carlo and quasi-Monte Carlo techniques for option pricing},
  in: Monte Carlo and quasi-Monte Carlo methods 1996
  (P. Hellekalek, G. Larcher, H. Niederreiter, and P. Zinterhof, eds.),
  Springer, Berlin, pp.~1--18 (1998)

\bibitem{AGHJLV16}
  %{A.~Ammon, T.~Hartung, K.~Jansen, H.~Le\"ovey, and J.~Vollmer},
  {Ammon, A., Hartung, T., Jansen, K., Le\"ovey, H., Vollmer, J.:}
  {On the efficient numerical solution of lattice systems with low-order couplings},
  {Comput.\ Phys.\ Commun.} %{Computer Physics Communications},
  {\bf 198}, 71--81 (2016)

\bibitem{CMO97}
  %{R.~E.~Caflisch, W.~Morokoff, and A.B.~Owen},
  {Caflisch, R.E., Morokoff, W., Owen, A.B.:}
  {Valuation of mortgage backed securities using Brownian bridges
  to reduce effective dimension},
  {J.~Comput.\ Finance} {\bf 1}, 27--46 (1997)

\bibitem{CD15}
  %{A.~Cohen and R.~DeVore},
  {Cohen, A., DeVore, R.:}
  {Approximation of high-dimensional parametric PDEs},
  {Acta Numer.} {\bf 24}, 1--159 (2015)

\bibitem{CKN06}
  %{R.~Cools, F.~Y.~Kuo, and D.~Nuyens},
  {Cools, R., Kuo, F.Y., Nuyens, D.:}
  {Constructing embedded lattice rules for multivariate integration},
  {SIAM J.~Sci.\ Comput.} {\bf 28}, 2162--2188 (2006)

\bibitem{CKNS16}
  %{R.~Cools, F.~Y.~Kuo, D.~Nuyens, and G.~Suryanarayana},
  {Cools, R., Kuo, F.Y., Nuyens, D., Suryanarayana, G.:}
  {Tent-transformed lattice rules for integration and approximation of
  multivariate non-periodic functions},
  {J.\ Complexity} {\bf 36}, 166--181 (2016)

\bibitem{DGLS}
  %{J.~Dick, R.~N.~Gantner, Q.~T.~Le~Gia, and Ch.~Schwab},
  {Dick, J., Gantner, R.N., Le Gia, Q.T., Schwab, Ch.:}
  {Higher order Quasi-Monte Carlo integration for Bayesian estimation}
  (in review)

\bibitem{DILP}
 %{J.~Dick, Ch.~Irrgeher, G.~Leobacher, and F.~Pillichshammer},
 {Dick, J., Irrgeher, Ch., Leobacher, G., Pillichshammer, F.:}
 {On the optimal order of integration in Hermite spaces with finite smoothness}
 (in review)

\bibitem{DKLNS14}
  %{J.~Dick, F.~Y.~Kuo, Q.~T.~Le~Gia, D.~Nuyens, and Ch.~Schwab},
  {Dick, J., Kuo, F.Y., Le Gia, Q.T., Nuyens, D., Schwab, Ch.:}
  {Higher order QMC Galerkin discretization for parametric
  operator equations},
  {SIAM J.\ Numer.\ Anal.} {\bf 52}, 2676--2702 (2014)

\bibitem{DKLS15}
  %{J.~Dick, F.~Y.~Kuo, Q.~T.~Le~Gia, and Ch.~Schwab},
  {Dick, J., Kuo, F.Y., Le Gia, Q.T., Schwab, Ch.:}
  {Fast QMC matrix-vector multiplication},
  {SIAM J.\ Sci.\ Comput.} {\bf 37}, A1436--A1450 (2015)

\bibitem{DKLS16}
  %{J.~Dick, F.~Y.~Kuo, Q.~T.~Le~Gia, and Ch.~Schwab},
  {Dick, J., Kuo, F.Y., Le Gia, Q.T., Schwab, Ch.:}
  {Multi-level higher order QMC Galerkin discretization for affine parametric
  operator equations},
  {SIAM J.\ Numer.\ Anal.}, {\bf 54}, 2541--2568 (2016)

\bibitem{DKS13}
  %{J.~Dick, F.~Y.~Kuo, and I.~H.~Sloan},
  {Dick, J., Kuo, F.Y., Sloan, I.H.:}
  {High-dimensional integration: the Quasi-Monte Carlo way},
  {Acta Numer.} \textbf{22}, 133--288 (2013)

\bibitem{DLS16}
  % {J.~Dick, Q.~T.~Le~Gia, and Ch.~Schwab},
  {Dick, J., Le Gia, Q.T., Schwab, Ch.:}
   Higher order Quasi-Monte Carlo integration for holomorphic, parametric operator equations,
   {SIAM/ASA Journal on Uncertainty Quantification}, {\bf 4}, 48--79 (2016)

\bibitem{DNP14}
 %J.~Dick, D.~Nuyens, and F.~Pillichshammer,
 {Dick, J., Nuyens, D., Pillichshammer, F.:}
 Lattice rules for nonperiodic smooth integrands,
 {Numer.\ Math.} \textbf{126}, 259--291 (2014)

\bibitem{DP10}
  %{J.~Dick and F.~Pillichshammer},
  {Dick, J., Pillichshammer, F.:}
  {Digital Nets and Sequences}, Cambridge University Press (2010)

\bibitem{DPW08}
  %{J.~Dick, F.~Pillichshammer, and B.~J.~Waterhouse},
  {Dick, J., Pillichshammer, F., Waterhouse, B.J.:}
  {The construction of good extensible rank-$1$ lattices},
  {Math.\ Comp.} {\bf 77}, 2345--2374 (2008)

\bibitem{DN97}
  %{C.~R.~Dietrich and G.~H.~Newsam},
  {Dietrich, C.R., Newsam, G.H.:}
  {Fast and exact simulation of stationary Gaussian processes through circulant
  embedding of the covariance matrix},
  {SIAM J.\ Sci.\ Comput.} \textbf{18}, 1088--1107 (1997)

\bibitem{FKS}
  %{M.~Feischl, F.~Y.~Kuo, and I.~H.~Sloan},
  {Feischl, M., Kuo, F.Y., Sloan, I.H.:}
  {Fast random field generation with $H$-matrices} (in review)

\bibitem{GH15}
  %{M.~Ganesh and S.~C.~Hawkins},
  {Ganesh, M., Hawkins, S.C.:}
  {A high performance computing and sensitivity analysis algorithm for stochastic
  many-particle wave scattering},
  {SIAM J.\ Sci.\ Comput.} {\bf 37}, A1475--A1503 (2015)

\bibitem{GanKS}
  %{M.~Ganesh, F.~Y.~Kuo, and I.~H.~Sloan},
  {Ganesh, M., Kuo, F.Y., Sloan, I.H.:}
  {Quasi-Monte Carlo finite element methods for stochastic heterogeneous wave progagation
  models}
  (in preparation)

\bibitem{GHS}
  %{R.~N.~Gantner, L.~Herrmann and Ch.~Schwab},
  {Gantner, R.N., Herrmann, L., Schwab, Ch.:}
  {Quasi-Monte Carlo integration for affine-parametric, elliptic PDEs: local supports imply product
  weights}
  (in review)

\bibitem{GGKSS}
  %{A.~D.~Gilbert, I.~G.~Graham, F.~Y.~Kuo, R.~Scheichl, and I.~H.~Sloan},
  {Gilbert, A.D., Graham, I.G., Kuo, F.Y., Scheichl, R., Sloan, I.H.:}
  {Quasi-Monte Carlo theory for an eigenproblem with a random coefficient}
  (in preparation)

\bibitem{GKNW}
  %{A.~D.~Gilbert, F.~Y.~Kuo, D.~Nuyens, and G.~W.~Wasilkowski},
  {Gilbert, A.D., Kuo, F.Y., Nuyens, D., Wasilkowski, G.W.:}
  {Efficient implementation of the multivariate decomposition method}
  (in preparation)

\bibitem{GilWas17}
  %{A.~D.~Gilbert and G.~W.~Wasilkowski},
  {Gilbert, A.D., Wasilkowski, G.W.:}
  {Small superposition dimension and active set construction for multivariate
  integration under modest error demand},
  {J.\ Complexity} {\bf 42}, 94--109 (2017)

\bibitem{Gil15}
  %{M.~B.~Giles},
  {Giles, M.B.:}
  {Multilevel Monte Carlo methods},
  {Acta Numer.} {\bf 24}, 259--328 (2015)

\bibitem{GilKS}
  %{M.~B.~Giles, F.~Y.~Kuo, and I.~H.~Sloan},
  {Giles, M.B., Kuo, F.Y., Sloan, I.H.:}
  {Combining sparse grids, multilevel MC and QMC for elliptic PDEs with random coefficients}
  (submitted to this volume)

\bibitem{GKSW08}
  %{M.~B.~Giles, F.~Y.~Kuo, I.~H.~Sloan, and B.~J.~Waterhouse},
  {Giles, M.B., Kuo, F.Y., Sloan, I.H., Waterhouse, B.J.:}
  {Quasi-Monte Carlo for finance applications},
  {ANZIAM J.} \textbf{50}, C308--C323 (CTAC2008) (2008)

\bibitem{GilWat09}
  %{M.~B.~Giles and B.~J.~Waterhouse},
  {Giles, M.B., Waterhouse, B.J.:}
  {Multilevel quasi-Monte Carlo path simulation}.
  {Radon Series Comp.\ Appl.\ Math.} {\bf 8}, 1--18 (2009)

\bibitem{God15}
  %{T.~Goda},
  {Goda, T.:}
  {Good interlaced polynomial lattice rules for
  numerical integration in weighted Walsh spaces},
  {J.\ Comput.\ Appl.\ Math.} {\bf 285}, 279--294 (2015)

\bibitem{GD15}
  %{T.~Goda and J.~Dick},
  {Goda, T., Dick, J.:}
  {Construction of interlaced scrambled polynomial lattice rules of arbitrary high order},
  {Found.\ Comput.\ Math.} {\bf 15}, 1245--1278 (2015)

\bibitem{GKNSSS15}
  %{I.~G.~Graham, F.~Y.~Kuo, J.~A.~Nichols, R.~Scheichl, Ch.~Schwab, and I.~H.~Sloan},
  {Graham, I.G., Kuo, F.Y., Nichols, J.A., Scheichl, R., Schwab, Ch., Sloan, I.H.:}
  {QMC FE methods for PDEs with log-normal random coefficients},
  {Numer.\ Math.} {\bf 131}, 329--368 (2015)

\bibitem{GKNSS11}
  %{I.~G.~Graham, F.~Y.~Kuo, D.~Nuyens, R.~Scheichl, and I.~H.~Sloan},
  {Graham, I.G., Kuo, F.Y., Nuyens, D., Scheichl, R., Sloan, I.H.:}
  {Quasi-Monte Carlo methods for elliptic PDEs with random coefficients and applications},
  {J.\ Comput.\ Phys.} {\bf 230}, 3668--3694 (2011)

\bibitem{GKNSS-paper1}
  %{I.~G.~Graham, F.~Y.~Kuo, D.~Nuyens, R.~Scheichl, and I.~H.~Sloan},
  {Graham, I.G., Kuo, F.Y., Nuyens, D., Scheichl, R., Sloan, I.H.:}
  {Analysis of circulant embedding methods for sampling stationary random fields}
  (in preparation)

\bibitem{GKNSS-paper2}
  %{I.~G.~Graham, F.~Y.~Kuo, D.~Nuyens, R.~Scheichl, and I.~H.~Sloan},
  {Graham, I.G., Kuo, F.Y., Nuyens, D., Scheichl, R., Sloan, I.H.:}
  {Circulant embedding with QMC -- analysis for elliptic PDE with lognormal coefficients}
  (in preparation)

\bibitem{GPS}
  %{I.~G.~Graham, M.~J.~Parkinson, and R.~Scheichl},
  {Graham, I.G., Parkinson, M.J., Scheichl, R.:}
  {Modern Monte Carlo variants for uncertainty quantification in neutron transport}
  (in review)

\bibitem{GKS13}
  %{M.~Griebel, F.~Y.~Kuo, and I.~H.~Sloan},
  {Griebel, M., Kuo, F.Y., Sloan, I.H.:}
  {The smoothing effect of integration in $\mathbb{R}^d$ and the ANOVA decomposition},
  {Math.\ Comp.} {\bf 82} (2013), 383--400;
  and the note in {Math.\ Comp.} {\bf 86}, 1847--1854 (2017)

\bibitem{GKLS}
  %{A.~Griewank, F.~Y.~Kuo, H.~Le\"ovey, and I.~H.~Sloan},
  {Griewank, A., Kuo, F.Y., Le\"ovey, H., Sloan, I.H.:}
  {High dimensional integration of kinks and jumps -- smoothing by preintegration}
  (in review)

\bibitem{GWZ14}
  %{M.~Gunzburger, C.~Webster, and G.~Zhang},
  {Gunzburger, M., Webster, C., Zhang, G.:}
  Stochastic finite element methods for partial differential equations with random input data,
  {Acta Numer.} {\bf 23}, 521--650 (2014)

\bibitem{Hac15}
  %{W.~Hackbusch},
  {Hackbusch, W:}
  {Hierarchical matrices: algorithms and analysis},
  Springer, Heidelberg (2015)

\bibitem{HNT16}
  %{A.~L.~Haji-Ali, F.~Nobile, and R.~Tempone},
  {Haji-Ali, A.L., Nobile, F., Tempone, R.:}
  {Multi-index Monte Carlo: when sparsity meets sampling},
  {Numer.\ Math.}  \textbf{132}, 767--806 (2016)

\bibitem{HPS17}
  %{H.~Harbrecht, M.~Peters, and M.~Siebenmorgen}.
  {Harbrecht, H., Peters, M., Siebenmorgen, M.:}
  {On the quasi-Monte Carlo method with Halton points for elliptic PDEs with log-normal
  diffusion},
  {Math.\ Comp.} {\bf 86}, 771--797 (2017)

\bibitem{HerSch}
  %{L.~Herrmann and Ch.~Schwab},
  {Herrmann, L., Schwab, Ch.:}
  {QMC integration for lognormal-parametric, elliptic PDEs: local supports imply product
  weights}
  (in review)

\bibitem{Hic02}
  %{F.~J.~Hickernell},
  {Hickernell, F.J.:}
  {Obtaining $O(N^{-2+\epsilon})$ convergence for lattice quadrature rules},
  in: Monte Carlo and Quasi-Monte Carlo Methods 2000
  (K.~T.~Fang, F.~J.~Hickernell, and H.~Niederreiter, eds.),
  Springer, Berlin, pp.~274--289 (2002)

\bibitem{HHLL00}
  %{F.~J.~Hickernell, H.~S.~Hong, P.~L\'Ecuyer, and C.~Lemieux},
  {Hickernell, F.J., Hong, H.S., L\'Ecuyer, P., Lemieux, C.:}
  {Extensible lattice sequences for quasi-Monte Carlo quadrature},
  {SIAM J.\ Sci.\ Comput.} \textbf{22}, 1117--1138 (2000)

\bibitem{HN03}
  %{F.~J.~Hickernell and H.~Niederreiter},
  {Hickernell, F.J., Niederreiter, H.:}
  {The existence of good extensible rank-$1$ lattices},
  {J. Complexity} \textbf{19}, 286--300 (2003)

\bibitem{JLGM14}
  %{K.~Jansen, H.~Le\"ovey, A.~Griewank, and M.~M\"uller-Preussker},
  {Jansen, K., Le\"ovey, H., Griewank, A., M\"uller-Preussker, M.:}
  {Quasi-Monte Carlo methods for lattice systems: a first look},
  {Comput.\ Phys.\ Commun.} {\bf 185}, 948--959 (2014)

\bibitem{Kaz}
  %{Y.~Kazashi},
  {Kazashi, Y.:}
  {Quasi-Monte Carlo integration with product weights for elliptic PDEs with log-normal
  coefficients}
  (in review)

\bibitem{KDSWW08}
  %{F.~Y.~Kuo, W.~D.~M.~Dunsmuir, I.~H.~Sloan, M.~P.~Wand, and R.~S.~Womersley},
  {Kuo, F.Y., Dunsmuir, W.D.M., Sloan, I.H., Wand, M.P., Womersley, R.S.:}
  {Quasi-Monte Carlo for highly sturctured generalised response models}
  %{Methodology and Computing in Applied Probability}
  {Methodol.\ Comput.\ App.} {\bf 10}, 239--275 (2008)

\bibitem{KN16}
  %{F.~Y.~Kuo and D.~Nuyens},
  {Kuo, F.Y., Nuyens, D.:}
  {Application of quasi-Monte Carlo methods to elliptic PDEs with
  random diffusion coefficients -- a survey of analysis and
  implementation},
  {Found.\ Comput.\ Math.} {\bf 16}, 1631--1696 (2016)

\bibitem{KN-tutorial}
  %{F.~Y.~Kuo and D.~Nuyens},
  {Kuo, F.Y., Nuyens, D.:}
  {Application of quasi-Monte Carlo methods to PDEs with random coefficients -- an overview and tutorial}
  (to appear in this volume)

\bibitem{KNPSW}
  %{F.~Y.~Kuo, D.~Nuyens, L.~Plaskota, I.~H.~Sloan, and G.~W.~Wasilkowski},
  {Kuo, F.Y., Nuyens, D., Plaskota, L., Sloan, I.H., Wasilkowski, G.W.:}
  {Infinite-dimensional integration and the multivariate decomposition
  method},
  J.\ Comput.\ Appl.\ Math. {\bf 326}, 217--234 (2017)

\bibitem{KSSSU}
  %{F.~Y.~Kuo, R.~Scheichl, Ch.~Schwab, I.~H.~Sloan, and E.~Ullmann},
  {Kuo, F.Y., Scheichl, R., Schwab, Ch., Sloan, I.H., Ullmann, E.:}
  {Multilevel Quasi-Monte Carlo methods for lognormal diffusion problems},
  {Math.\ Comp.} {\bf 86}, 2827--2860 (2017)

\bibitem{KSS11}
  %{F.~Y.~Kuo, Ch.~Schwab, and I.~H.~Sloan},
  {Kuo, F.Y., Schwab, Ch., Sloan, I.H.:}
  {Quasi-Monte Carlo methods for high dimensional integration:
  the standard weighted-space setting and beyond},
  {ANZIAM J.} {\bf 53}, 1--37 (2011)

\bibitem{KSS12}
  %{F.~Y.~Kuo, Ch.~Schwab, and I.~H.~Sloan},
  {Kuo, F.Y., Schwab, Ch., Sloan, I.H.:}
  {Quasi-Monte Carlo finite element methods for a class of elliptic partial
  differential equations with random coefficient},
  {SIAM J.\ Numer.\ Anal.} {\bf 50}, 3351--3374 (2012)

\bibitem{KSS15}
  %{F.~Y.~Kuo, Ch.~Schwab, and I.~H.~Sloan},
  {Kuo, F.Y., Schwab, Ch., Sloan, I.H.:}
  {Multi-level quasi-Monte Carlo finite element methods for a class of elliptic partial
  differential equations with random coefficient},
  {Found.\ Comput.\ Math.} {\bf 15}, 411--449 (2015)

\bibitem{KSWWat10}
  %{F.~Y.~Kuo, I.~H.~Sloan, G.~W.~Wasilkowski, and B.~J.~Waterhouse},
  {Kuo, F.Y., Sloan, I.H., Wasilkowski, G.W., Waterhouse, B.J.:}
  {Randomly shifted lattice rules with the optimal rate of convergence for unbounded integrands},
  {J.\ Complexity} {\bf 26}, 135--160 (2010)

\bibitem{KSWW10a}
  %{F.~Y.~Kuo, I.~H.~Sloan, G.~W.~Wasilkowski, and H.~Wo\'zniakowski},
  {Kuo, F.Y., Sloan, I.H., Wasilkowski, G.W., Wo\'zniakowski, H.:}
  {Liberating the dimension},
  {J.\ Complexity} {\bf 26}, 422--454 (2010)

\bibitem{KSWW10b}
  %{F.~Y.~Kuo, I.~H.~Sloan, G.~W.~Wasilkowski, and H.~Wo\'zniakowski},
  {Kuo, F.Y., Sloan, I.H., Wasilkowski, G.W., Wo\'zniakowski, H.:}
  {On decompositions of multivariate functions},
  {Math.\ Comp.} {\bf 79}, 953--966 (2010)

\bibitem{Lem09}
  %{C.~Lemieux},
  {Lemieux, C.:}
  {Monte Carlo and Quasi-Monte Carlo Sampling},
  Springer, New York (2009)

\bibitem{LP14}
  %{G.~Leobacher and F.~Pillichshammer},
  {Leobacher, G., Pillichshammer, F.:}
  {Introduction to Quasi-Monte Carlo Integration and Applications},
  Springer (2014)

\bibitem{NN}
  %{D.~T.~P.~Nguyen and D.~Nuyens},
  {Nguyen, D.T.P., Nuyens, D.:}
  {Multivariate integration over $\mathbb{R}^s$ with exponential rate of convergence},
  {J.\ Comput.\ Appl.\ Math.} {\bf315}, 327--342 (2017)

\bibitem{NK14}
  %{J.~A.~Nichols and F.~Y.~Kuo},
  {Nichols, J.A., Kuo F.Y.:}
  {Fast CBC construction of randomly shifted lattice rules
  achieving $\mathcal{O}(N^{-1+\delta})$ convergence for unbounded integrands
  in $\mathbb{R}^s$ in weighted spaces with POD weights},
  {J.~Complexity} {\bf 30}, 444--468 (2014)

\bibitem{Nie92}
  %{H.~Niederreiter},
  {Niederreiter, H.:}
  {Random Number Generation and Quasi-Monte Carlo Methods},
  SIAM, Philadelphia (1992)

\bibitem{NW10}
  %{E.~Novak and H.~Wo\'zniakowski},
  {Novak, E., Wo\'zniakowski, H.:}
  {Tractability of Multivariate Problems, II: Standard Information for Functionals},
  European Mathematical Society, Z\"urich (2010)

\bibitem{Nuy14}
  %{D.~Nuyens},
  {Nuyens, D.:}
  The construction of good lattice rules and polynomial lattice rules,
  in: Uniform Distribution and Quasi-{M}onte {C}arlo Methods
  (P.~Kritzer, H.~Niederreiter, F.~Pillichshammer, A.~Winterhof, eds.),
  Radon Series on Computational and Applied Mathematics Vol.~15,
  De Gruyter, pp.~223--256 (2014)

\bibitem{NC06a}
  %{D.~Nuyens and R.~Cools},
  {Nuyens, D., Cools, R.:}
  {Fast algorithms for component-by-component construction of
  rank-$1$ lattice rules in shift-invariant reproducing kernel
  Hilbert spaces},
  {Math.\ Comp.} {\bf 75}, 903--920 (2006)

\bibitem{NC06b}
  %{D.~Nuyens and R.~Cools},
  {Nuyens, D., Cools, R.:}
  {Fast component-by-component construction of rank-$1$ lattice rules with
  a non-prime number of points},
  {J.\ Complexity} {\bf 22}, 4--28 (2006)

\bibitem{Owe97}
  %{A.~B.~Owen},
  {Owen, A.B.:}
  {Scrambled net variance for integrals of smooth functions},
  {Ann.\ Statist.} \textbf{25}, 1541--1562 (1997)

\bibitem{Owe06}
  %{A.~B.~Owen},
  {Owen, A.B.:}
  {Halton sequences avoid the origin},
  {SIAM Rev.} \textbf{48}, 487--503 (2006)

\bibitem{RNV}
  %{P.~Robbe, D.~Nuyens, and S.~Vandewalle},
  {Robbe, P., Nuyens, D., Vandewalle, S.:}
  {A multi-index quasi-Monte Carlo algorithm for lognormal diffusion
  problems},
  SIAM J.\ Sci.\ Comput.\ (to appear)

\bibitem{SST17}
  %{R.~Scheichl, A.~Stuart, and A.~L.~Teckentrup},
  {Scheichl, R., Stuart, A., Teckentrup, A.L.:}
  {Quasi-Monte Carlo and multilevel Monte Carlo methods for computing posterior
  expectations in elliptic inverse problems},
  {SIAM/ASA Journal on Uncertainty Quantification}
  {\bf 5}, 493--518 (2017)

\bibitem{SG11}
  %{Ch.~Schwab and C.~J.~Gittelson},
  {Schwab, Ch., Gittelson, C.J.:}
  {Sparse tensor discretizations of high-dimensional parametric and stoch
  astic PDEs},
  {Acta Numer.} {\bf 20}, 291--467 (2011)

\bibitem{SKS13}
  %{V.~Sinescu, F.~Y.~Kuo, and I.~H. Sloan},
  {Sinescu, V., Kuo, F.Y., Sloan, I.H.:}
  {On the choice of weights in a function space for quasi-Monte Carlo
  methods for a class of generalised response models in statistics},
  in: Monte Carlo and Quasi-Monte Carlo Methods 2012
  (J.~Dick, F.~Y.~Kuo, G.~Peters, and I.~H.~Sloan, eds),
  Springer Verlag, Heidelberg, pp.~631--647 (2013)

\bibitem{SJ94}
  %{I.~H.~Sloan and S.~Joe},
  {Sloan, I.H., Joe, S.:}
  {Lattice Methods for Multiple Integration},
  Oxford University Press, Oxford (1994)

\bibitem{SKJ02b}
  %{I.~H.~Sloan, F.~Y.~Kuo, and S.~Joe},
  {Sloan, I.H., Kuo, F.Y., Joe, S.:}
  {Constructing randomly shifted lattice rules in weighted Sobolev
  spaces},
  {SIAM J.~Numer.\ Anal.} {\bf 40}, 1650--1665 (2002)

\bibitem{SWW04}
  %{I.~H.~Sloan, X.~Wang, and H.~Wo\'zniakowski},
  {Sloan, I.H., Wang, X., Wo\'zniakowski, H.:}
  {Finite-order weights imply tractability of multivariate integration},
  {J.~Complexity} {\bf 20}, 46--74 (2004)

\bibitem{SW98}
  %{I.~H.~Sloan and H.~Wo\'zniakowski},
  {Sloan, I.H., Wo\'zniakowski, H.:}
  {When are quasi-Monte Carlo algorithms efficient for
  high-dimensional integrals?},
  {J.~Complexity} {\bf 14}, 1--33 (1998)

\bibitem{Was13}
  %{G.~W.~Wasilkowski},
  {Wasilkowski, G.W.:}
  {On tractability of linear tensor product problems for $\infty$-variate classes of functions},
  {J.\ Complexity}, {\bf 29}, 351--369 (2013)
\end{thebibliography}
% and then copy paste the contents of the .bbl file here for the final version.
%
% E.g.:
%\begin{thebibliography}{99.}%

%\bibitem{ACN2013}
%N.~Achtsis, R.~Cools and D.~Nuyens.
%\newblock Conditional sampling for barrier option pricing under the Heston model.
%\newblock In J.~Dick, F.~Y.\ Kuo, G.~W.\ Peters and I.~H.\ Sloan, editors, {\em {M}onte {C}arlo
%  and Quasi-{M}onte {C}arlo Methods 2012}, pages 253--269. Springer-Verlag, 2013.

%\bibitem{CKN2006}
%R.~Cools, F.~Y. Kuo, and D.~Nuyens.
%\newblock Constructing embedded lattice rules for multivariate integration.
%\newblock {\em SIAM Journal on Scientific Computing}, 28(6):2162--2188, 2006.

%\bibitem{DP2010}
%J.~Dick and F.~Pillichshammer.
%\newblock {\em Digital Nets and Sequences: Discrepancy Theory and Quasi-Monte
%  Carlo Integration}.
%\newblock Cambridge University Press, 2010.

%\bibitem{IT2006}
%J.~{Imai} and K.~S.\ {Tan}.
%\newblock A general dimension reduction technique for derivative pricing.
%\newblock {\em Journal of Computational Finance}, 10(2):129--155, 2006.

%\bibitem{LEC2009}
% P.~L'{\'E}cuyer.
%\newblock Quasi-Monte Carlo methods with applications in finance.
%\newblock {\em Finance and Stochastics}, 13(3):307--349, 2009.

%\end{thebibliography}

\end{document}